\documentclass{amsart}
\usepackage[utf8]{inputenc}
\usepackage{zeyu}
\usepackage{calc}
\title{The Oddtown problem modulo a composite number}

\makeatletter
\def\paragraph{\@startsection{paragraph}{4}%
  \z@\z@{-\fontdimen2\font}%
  {\normalfont\bfseries}}
\makeatother

\newcommand*{\spn}[1]{\operatorname{span}_{#1}}                 
\DeclareMathOperator{\PG}{PG}
\DeclareMathOperator{\diag}{diag}

\begin{document}
\begin{abstract}
A family of subsets $\cA$ of an $n$-element set is called an $\ell$-Oddtown if the sizes of all sets are not divisible by $\ell$, but the sizes of pairwise intersections are divisible by~$\ell$.
Berlekamp and Graver showed that when $\ell$ is a prime, the maximum size of an $\ell$-Oddtown is~$n$.
Babai and Frankl extended this to prime powers, and asked whether the maximum size is still $n$ when $\ell$ is not a prime power, a question that was open even for $\ell=6$.
For square-free composite moduli with $\omega$ distinct prime factors, the argument of Szegedy gives an upper bound of
$\omega n-\omega\log_2 n$ on the size of an $\ell$-Oddtown.
We answer the question of Babai and Frankl in the negative by constructing $\ell$-Oddtowns of size $\omega n-o(n)$, which shows that the leading term $\omega n$ cannot be improved.
We also improve Szegedy's upper bound to
$\omega n-(2\omega +\varepsilon)\log_2 n$ for most $\ell$ and $n$ using a combination of linear algebraic and Fourier-analytic arguments.

\end{abstract}
\maketitle

\section{Introduction}

The Oddtown problem is perhaps the simplest application of the linear algebra method in combinatorics,
and it is often the first use of the method that students encounter. The problem asks for the largest
size of a family $\cA$ of subsets of $[n]:= \{1,2,\dotsc,n\}$ of odd cardinality such that the intersection of any two sets in the family is of even size.
The original motivation for the question was a problem of Erd\H{o}s that was solved independently by Berlekamp \cite{berlekamp1969subsets} and Graver \cite{graver1975boolean}.
The name of the problem comes from the book of Babai and Frankl \cite{babai1992linear}, where it is cast
as a story about a city whose council tries to contain proliferation of clubs formed by its eccentric residents.

After identifying each set with its characteristic vector, simple mod-$2$ linear algebra shows that
the family $\cA$ can have no more than $n$ sets. The bound is sharp, as witnessed by the family of all single-element
sets. There are many other sharp examples \cite[Ex~1.1.14]{babai1992linear}.

There have been a number of works examining applications of linear algebra 
modulo a composite number to combinatorics. A particularly notable one is the line of works by Grolmusz \cite{grolmusz2000superpolynomial,grolmusz2003note}
who constructed large families of sets with nice intersection properties modulo a composite number.
These found numerous applications, and our construction in \Cref{sec:lower} is in the same spirit.
Another adjacent recent line of works \cite{arsovski2021p,dhar2020proof,dhar2022maximal} concerns Kakeya problem over $\ZZ/{N\ZZ}$,
which is a toy model for the original Kakeya problem.

Motivated by the desire to better understand linear algebra arguments in combinatorics, in this paper
we examine the generalization of the Oddtown problem to composite moduli. By \emph{$\ell$-Oddtown} we mean
a set family $\cA\subset 2^{[n]}$ satisfying
\begin{itemize}
  \item $\abs{A}\nequiv 0\pmod \ell$ for any $A\in \cA$, and
  \item $\abs{A\cap B}\equiv 0\pmod \ell$ for any distinct $A,B\in \cA$. 
\end{itemize}

Let $f_{\ell}(n)$ denote the maximum size of an $\ell$-Oddtown.

\paragraph{Known bounds.} Let $\omega$ be the number of distinct prime factors
of $\ell$. If $\omega=1$, i.e., $\ell$ is a power of a prime,
the same sharp bound $f_{\ell}(n)\leq n$ can be proved using basic linear algebra \cite{babai1980set}.

More generally, let $\ell=p_1^{\alpha_1}\dotsb p_{\omega}^{\alpha_{\omega}}$ be the prime factorization
of $\ell$, and set $\cA_i:= \{A\in \cA : \abs{A}\nequiv 0\pmod {p_i^{\alpha_i}}\}$; then
$\abs{\cA_i}\leq n$ by the prime-power modulus bound just mentioned, which in view of $\cA=\cA_1\cup\dotsb \cup \cA_{\omega}$ implies what we call the \emph{trivial bound}
\[
  f_{\ell}(n)\leq \omega n.
\]
As already mentioned, for $\omega=1$ this bound is sharp. So, for the rest of the paper we assume that $\omega\geq 2$.

Babai and Frankl \cite[Ex~1.1.27]{babai1992linear} asked whether the sharp bound $f_{\ell}(n)\leq n$
persists when $\ell$ is not a prime power, and noted that this question was open even for $\ell=6$.
The question was reiterated by Sudakov and Vieira \cite{sudakov2016two}.

The only improvement to the trivial bound that we are aware of is due to Mario Szegedy \cite[Ex~1.1.28]{babai1992linear}.
The statement appearing in \cite{babai1992linear} is for $\ell=6$ only, but the suggested argument generalizes to
\begin{equation}\label{eq:szegedy}
  f_{\ell}(n)\leq \omega n -\omega \log n\qquad\text{for square-free }\ell,
\end{equation}
Here and throughout the paper, $\log$ denotes the base-$2$ logarithm.

\paragraph{New results.} Our first result answers the question of Babai and Frankl, and shows that the
constant $\omega$ in the trivial bound cannot be improved.
\begin{theorem}\label{thm:construction}
  Let $\ell$ be a positive integer with at least two distinct prime factors. If $n$ is sufficiently
  large compared to $\ell$, then
  \begin{equation*}
    f_{\ell}(n)\geq \omega n-O_{\ell}\Bigl(n^{\frac{\omega-2}{\omega-1}}(\log n)^{C}\Bigr),
  \end{equation*}
  where $C=C(\ell)$. In particular, $f_{\ell}(n)=(\omega+o(1))n$.
\end{theorem}

The construction has two ingredients. First, for each $i$ we exhibit an $\ell$-Oddtown of size $n$ such
that the matrix formed by the characteristic vectors of its sets has sublinear inner rank over
$\ZZ/p_t^{\alpha_t}\ZZ$ for every $t\neq i$, built from complements of hyperplanes in projective
spaces over fields of characteristic $p_t$.
Second, this low rank in the wrong characteristics allows us to superimpose the $\omega$ Oddtowns so
obtained on a single ground set at the cost of few additional elements, by a gadget based on the
Chinese remainder theorem.

In the other direction, we more than double the $\omega \log n$ saving in \eqref{eq:szegedy}, and
extend the result to general (non-square-free) moduli.
\begin{theorem}\label{thm:main}
  Let $\ell$ be a positive integer with at least two distinct prime factors. Then for any $n\geq 4$,
  \begin{equation*}
    f_{\ell}(n)\leq \omega n - 2 \omega\log n+11.
  \end{equation*}
\end{theorem}

\begin{theorem}\label{thm:main2}
  Let $\ell$ be a positive integer with at least two distinct odd prime factors. If $n$ is sufficiently large compared to $\ell$, then 
  \begin{equation*}
    f_{\ell}(n)\leq \omega n - (2\omega+\varepsilon) \log n,
  \end{equation*}
  where $\varepsilon=\frac{1}{20}(\sum_{p\in P_\text{odd}^*(\ell)} p^{-2})$, and $P_{\text{odd}}^*(\ell)$ is the set of all
  odd primes $p$ dividing  $\ell$ \emph{except} for the smallest such~prime. 
\end{theorem}

The upper bounds are a combination of two new ingredients. First, we observe that in Szegedy's argument certain pairs of subspaces
are orthogonal to each other. This permits us to replace $\omega \log n$ by $2\omega \log n$, but only for square-free moduli.
The main difficulty is to generalize the orthogonality argument from prime moduli to prime-power moduli. This is done in \Cref{lem:dimension}.

\paragraph{Ranks of $0$--$1$ matrices with distinct rows and distinct columns.}
Much of the paper is spent improving this bound by further $\varepsilon \log n$. In the proof of \Cref{eq:szegedy}, a key
role is played by the bound $\abs{\{0,1\}^n\cap W}\leq 2^d$ that holds for any vector subspace $W$ of $\FF_p^n$ of dimension~$d$.
Equivalently, if $M\in\FF_p^{r\times c}$ is a $0$--$1$ matrix with 
distinct rows, then $d:=\rank_{p}M$ is always lower bounded by 
$\log r$.

It follows that, if \emph{both} rows and columns of $M$ are distinct,
then $d\geq \max(\log r,\log c)$. The main ingredient behind the $\varepsilon \log n$
improvement is the following result that excludes the possibility that both $\log r$ and $\log c$ are close to $d$ in this case.
\begin{mainlemma}[Informal version]
Let $p$ be an odd prime and let $M\in \FF_p^{r\times c}$ be a $0$--$1$ matrix with distinct columns and distinct rows. Set $d:= \rank_p M$. Then we have
\[
  d\gtrsim \min\left(\left(1+\frac{1}{36p^2}\right)\log r,3\log c\right).
\]
\end{mainlemma}
The improvement is depicted in the following figure.

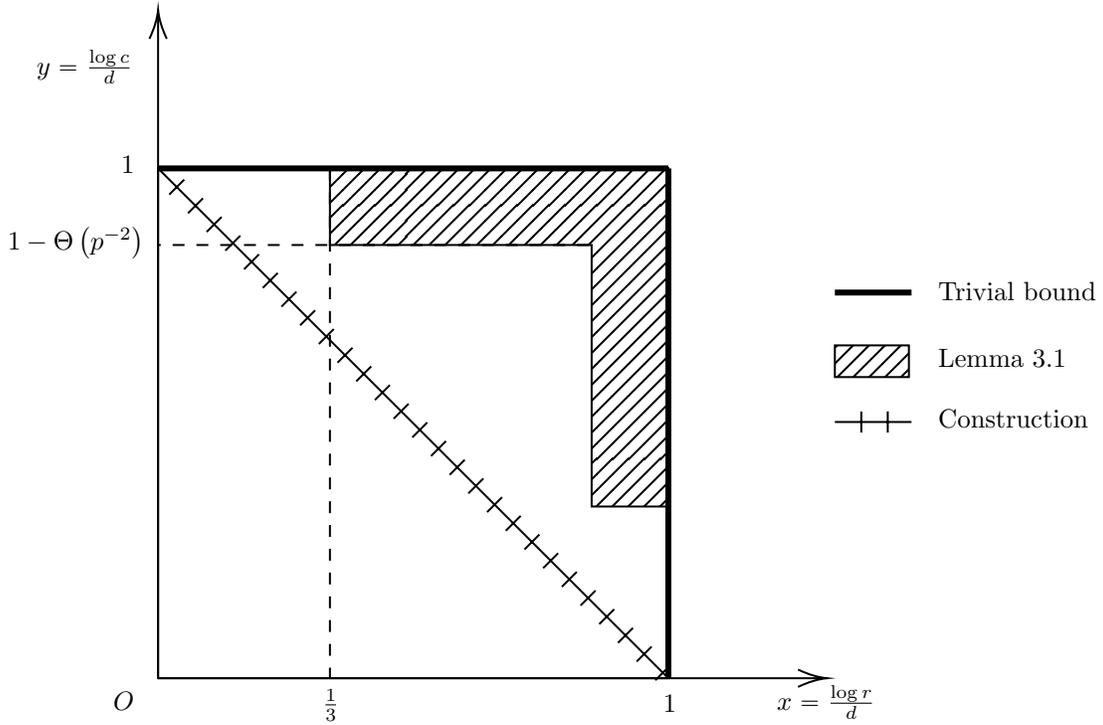
\begin{figure}[H]

 
\tikzset{
pattern size/.store in=\mcSize, 
pattern size = 5pt,
pattern thickness/.store in=\mcThickness, 
pattern thickness = 0.3pt,
pattern radius/.store in=\mcRadius, 
pattern radius = 1pt}
\makeatletter
\pgfutil@ifundefined{pgf@pattern@name@_0mzg5sxqy}{
\pgfdeclarepatternformonly[\mcThickness,\mcSize]{_0mzg5sxqy}
{\pgfqpoint{0pt}{0pt}}
{\pgfpoint{\mcSize+\mcThickness}{\mcSize+\mcThickness}}
{\pgfpoint{\mcSize}{\mcSize}}
{
\pgfsetcolor{\tikz@pattern@color}
\pgfsetlinewidth{\mcThickness}
\pgfpathmoveto{\pgfqpoint{0pt}{0pt}}
\pgfpathlineto{\pgfpoint{\mcSize+\mcThickness}{\mcSize+\mcThickness}}
\pgfusepath{stroke}
}}
\makeatother

 
\tikzset{
pattern size/.store in=\mcSize, 
pattern size = 5pt,
pattern thickness/.store in=\mcThickness, 
pattern thickness = 0.3pt,
pattern radius/.store in=\mcRadius, 
pattern radius = 1pt}
\makeatletter
\pgfutil@ifundefined{pgf@pattern@name@_5ge3dieaw}{
\pgfdeclarepatternformonly[\mcThickness,\mcSize]{_5ge3dieaw}
{\pgfqpoint{0pt}{0pt}}
{\pgfpoint{\mcSize+\mcThickness}{\mcSize+\mcThickness}}
{\pgfpoint{\mcSize}{\mcSize}}
{
\pgfsetcolor{\tikz@pattern@color}
\pgfsetlinewidth{\mcThickness}
\pgfpathmoveto{\pgfqpoint{0pt}{0pt}}
\pgfpathlineto{\pgfpoint{\mcSize+\mcThickness}{\mcSize+\mcThickness}}
\pgfusepath{stroke}
}}
\makeatother
\tikzset{every picture/.style={line width=0.75pt}} 

\begin{tikzpicture}[x=1pt,y=1pt,yscale=-1,xscale=1]

\draw  [dash pattern={on 4.5pt off 4.5pt}]  (172,91) -- (172,284) ;
\draw   (107,91) -- (300,91) -- (300,284) -- (107,284) -- cycle ;
\draw    (107,91) -- (300,284) (116.9,95.24) -- (111.24,100.9)(123.97,102.31) -- (118.31,107.97)(131.04,109.38) -- (125.38,115.04)(138.11,116.46) -- (132.46,122.11)(145.18,123.53) -- (139.53,129.18)(152.25,130.6) -- (146.6,136.25)(159.33,137.67) -- (153.67,143.33)(166.4,144.74) -- (160.74,150.4)(173.47,151.81) -- (167.81,157.47)(180.54,158.88) -- (174.88,164.54)(187.61,165.95) -- (181.95,171.61)(194.68,173.02) -- (189.02,178.68)(201.75,180.1) -- (196.1,185.75)(208.82,187.17) -- (203.17,192.82)(215.89,194.24) -- (210.24,199.89)(222.97,201.31) -- (217.31,206.97)(230.04,208.38) -- (224.38,214.04)(237.11,215.45) -- (231.45,221.11)(244.18,222.52) -- (238.52,228.18)(251.25,229.59) -- (245.59,235.25)(258.32,236.66) -- (252.66,242.32)(265.39,243.74) -- (259.74,249.39)(272.46,250.81) -- (266.81,256.46)(279.53,257.88) -- (273.88,263.53)(286.61,264.95) -- (280.95,270.61)(293.68,272.02) -- (288.02,277.68)(300.75,279.09) -- (295.09,284.75) ;
\draw    (107,284) -- (107,33) ;
\draw [shift={(107,31)}, rotate = 90] [color={rgb, 255:red, 0; green, 0; blue, 0 }  ][line width=0.75]    (10.93,-3.29) .. controls (6.95,-1.4) and (3.31,-0.3) .. (0,0) .. controls (3.31,0.3) and (6.95,1.4) .. (10.93,3.29)   ;
\draw    (107,284) -- (358,284) ;
\draw [shift={(360,284)}, rotate = 180] [color={rgb, 255:red, 0; green, 0; blue, 0 }  ][line width=0.75]    (10.93,-3.29) .. controls (6.95,-1.4) and (3.31,-0.3) .. (0,0) .. controls (3.31,0.3) and (6.95,1.4) .. (10.93,3.29)   ;
\draw  [pattern=_0mzg5sxqy,pattern size=6pt,pattern thickness=0.75pt,pattern radius=0pt, pattern color={rgb, 255:red, 0; green, 0; blue, 0}] (300,219) -- (271,219) -- (271,120) -- (172,120) -- (172,91) -- (300,91) -- cycle ;
\draw  [dash pattern={on 4.5pt off 4.5pt}]  (270,120) -- (107,120) ;
\draw [line width=2.25]    (107,91) -- (300,91) ;
\draw [line width=2.25]    (300,284) -- (300,91) ;
\draw  [pattern=_5ge3dieaw,pattern size=6pt,pattern thickness=0.75pt,pattern radius=0pt, pattern color={rgb, 255:red, 0; green, 0; blue, 0}] (363,158) -- (391,158) -- (391,170) -- (363,170) -- cycle ;
\draw [line width=2.25]    (363,138) -- (392,138) ;
\draw    (363,187) -- (392,187) (373,183) -- (373,191)(383,183) -- (383,191) ;

\draw (60,45) node [anchor=north west][inner sep=0.75pt]  [font=\normalsize]  {$y=\frac{\log c}{d}$};
\draw (340,285.5) node [anchor=north west][inner sep=0.75pt]  [font=\normalsize]  {$x=\frac{\log r}{d}$};
\draw (168,287) node [anchor=north west][inner sep=0.75pt]  [font=\normalsize]  {$\frac{1}{3}$};
\draw (49,112) node [anchor=north west][inner sep=0.75pt]  [font=\normalsize]  {$1-\Theta \left( p^{-2}\right)$};
\draw (297,289) node [anchor=north west][inner sep=0.75pt]  [font=\normalsize]  {$1$};
\draw (89,288) node [anchor=north west][inner sep=0.75pt]  [font=\normalsize]  {$O$};
\draw (401,158) node [anchor=north west][inner sep=0.75pt]   [align=left] {\Cref{lem:01matrix}};
\draw (401,133) node [anchor=north west][inner sep=0.75pt]   [align=left] {Trivial bound};
\draw (401,181) node [anchor=north west][inner sep=0.75pt]   [align=left] {Construction};
\draw (92,85) node [anchor=north west][inner sep=0.75pt]  [font=\normalsize]  {$1$};

\end{tikzpicture}
\end{figure}
It is not difficult to construct $r$-by-$c$ matrices with distinct rows
and distinct columns of rank at most $\approx \log r+\log c$. Indeed, let $M_n\in \FF_p^{n\times 2^n}$ 
be the $0$--$1$ matrix whose columns are all possible $0$--$1$~vectors of length~$n$.
Then $M=\left(\begin{smallmatrix}M_{a}^{\top}&0\\0&M_b\end{smallmatrix}\right)$
is an $(2^a+b)$-by-$(a+2^b)$ matrix. Since $\rank_p M_n= n$, the matrix $M$ is of rank $a+b$. 

We suspect that the fraction $1/36p^2$ in \Cref{lem:01matrix} can be improved to  some absolute constant $\varepsilon'>0$.
This would imply that the constant $\varepsilon>0$ in \Cref{thm:main2} can be made independent
of~$\ell$.

\paragraph{Paper organization.}
The remainder of this paper is structured as follows. In \Cref{sec:lower} we give the construction proving \Cref{thm:construction}.
In \Cref{sec:algebra} we introduce our algebraic tools, and prove
\Cref{thm:main}. In \Cref{sec:combinatorics} we present the proof of \Cref{thm:main2}, modulo the proof of \Cref{lem:01matrix}, which is postponed to \Cref{sec:analysis}.
\Cref{lem:01matrix} is a combination of the Fourier-analytic argument in \Cref{pro:fourier} and a combinatorial covering argument in \Cref{lem:subMatrix}. Unlike the combinatorial part
of the argument, the Fourier-analytic component is rather unsophisticated. However, it is interesting because the applications of Fourier analysis to extremal set theory are very rare.
The most recent example is probably an approach to the density Hales--Jewett theorem by Gowers and Karam \cite{gowers2023low}.

\paragraph{Disclosure of AI use.}
The lower-bound construction in \Cref{sec:lower} was first proposed by GPT-5.6 Sol in response to prompts from the authors. The upper-bound results were obtained without AI assistance.

\section{The lower bound}\label{sec:lower}

We identify a set with its characteristic vector, and collect the characteristic vectors of a family
$\cA=\{A_1,\dotsc,A_m\}$ of subsets of $[n]$ as the rows of an $m$-by-$n$ $0$--$1$ matrix $M$. Since
$(MM^{\top})_{ij}=\abs{A_i\cap A_j}$, the family $\cA$ is an $\ell$-Oddtown precisely when
\[ (MM^{\top})_{ii}\nequiv0\pmod{\ell}\ \text{ for all } i, \qquad (MM^{\top})_{ij}\equiv0\pmod{\ell}\ \text{ for all } i\neq j.
\]
Adjoining a new ground element amounts to appending a column to $M$.

For each prime-power factor $p_i^{\alpha_i}$ of $\ell$, we construct a square $0$--$1$ matrix $M_i$
whose rows form an $\ell$-Oddtown and whose inner rank is small modulo every other prime-power factor
of $\ell$. Complements of projective hyperplanes provide the basic matrices, tensor products impose the
required congruences, and block-diagonal sums give matrices of every sufficiently large order. We then
place $M_1,\dotsc,M_{\omega}$ on the same columns. Their mixed Gram blocks have low rank, so a small
number of additional columns can cancel them modulo $\ell$. This superimposition idea goes back to
Grolmusz \cite{grolmusz2000superpolynomial}.

\subsection{Inner rank}\label{sec:rank}
The number of columns in the final correction is controlled by factorization ranks over the rings
$\ZZ/p^{\alpha}\ZZ$. We begin by recording the relevant notion of rank and the properties used below.

\begin{definition}\label{def:rank}
  Let $R$ be a commutative ring and let $M$ be a matrix over $R$. We write $\rank_R(M)$ for the least
  $d$ such that $M=XY$, where $X$ has $d$ columns and $Y$ has $d$ rows. This is the \emph{inner rank}
  of $M$ in the terminology of Cohn \cite[p.~3]{cohn2006free}.
\end{definition}

Over a field the inner rank is the usual rank, so that $\rank_{\FF_p}(M)=\rank_p M$ in the notation
used elsewhere in this paper. In general, $\rank_R(M)$ is at most the number of rows of $M$ and at most
the number of its columns, since $M=IM=MI$ for identity matrices of the appropriate sizes.

\begin{remark}\label{rem:whyRingRank}
  The factorization must take place over $\ZZ/p^{\alpha}\ZZ$, rather than merely over $\FF_p$. For
  example, $pI_m$ vanishes modulo $p$, whereas
  $\rank_{\ZZ/p^2\ZZ}(pI_m)\geq m/2$. Indeed, if $pI_m=XY$ with inner dimension $d$, then the columns
  of $X$ generate a submodule containing $p(\ZZ/p^2\ZZ)^m$. The former has at most $p^{2d}$ elements,
  while the latter has $p^m$, so $2d\geq m$. Thus reduction modulo $p$ can hide large factorization
  rank over $\ZZ/p^{\alpha}\ZZ$, which is why the estimates below must be carried out over the full
  prime-power ring.
\end{remark}

We shall use the following properties, all immediate from \Cref{def:rank}.

\begin{proposition}\label{pro:rankProperties}
  Let $R$ be a commutative ring, and write $\otimes$ for the Kronecker product.
  \begin{enumerate}[(i)]
    \item $\rank_R(A^{\top})=\rank_R(A)$.
    \item $\rank_R(AB)\leq\min(\rank_R(A),\rank_R(B))$.
    \item $\rank_R(A\otimes B)\leq \rank_R(A)\rank_R(B)$.
    \item If $A=\diag(A_1,\dotsc,A_k)$ is block diagonal, then
      $\rank_R(A)\leq\sum_{i=1}^{k}\rank_R(A_i)$.
  \end{enumerate}
\end{proposition}

\begin{proof}
  Part (i) follows by transposing a factorization. For (ii), a factorization $A=XY$ gives
  $AB=X(YB)$, and a factorization $B=XY$ gives $AB=(AX)Y$. For (iii), factorizations $A=XY$ and
  $B=UV$ give $A\otimes B=(X\otimes U)(Y\otimes V)$. For (iv), place minimum-width factorizations of
  the diagonal blocks block diagonally.
\end{proof}

\subsection{Complements of hyperplanes}\label{sec:hyperplanes}
Let $q$ be a prime power and let $s\geq2$. Write $\bar{x}$ for the point of $\PG(s-1,q)$ with
homogeneous coordinates $x\in\FF_q^s\setminus\{0\}$, and equip $\FF_q^s$ with the standard bilinear
form
\[
  x\cdot y=x_1y_1+\dotsb+x_sy_s.
\]
The number of points of $\PG(s-1,q)$ is
\[
  \frac{q^{s}-1}{q-1}=1+q+\dotsb+q^{s-1}\leq q^s.
\]
Let $N(s,q)$ be the $0$--$1$ matrix whose rows and columns are indexed by the points of
$\PG(s-1,q)$ and whose entry in position $(\bar{x},\bar{y})$ is $1$ if and only if $x\cdot y\neq0$.
This is well defined because rescaling $x$ or $y$ by a nonzero scalar does not affect whether
$x\cdot y$ vanishes. The row indexed by $\bar{x}$ is the characteristic vector of the complement of
the hyperplane $\{\bar{y}:x\cdot y=0\}$. These matrices have the two properties needed later: their row
weights and pairwise dot products are explicit, which allows us to control their divisibility, while
their inner rank over $\ZZ/p^{\alpha}\ZZ$ grows only polynomially in $s$ when $q$ is a power of $p$.

\begin{lemma}\label{lem:projectiveCounts}
  Every row of $N(s,q)$ has exactly $q^{s-1}$ ones, and any two distinct rows have exactly
  $(q-1)q^{s-2}$ ones in common.
\end{lemma}

\begin{proof}
  Fix a nonzero $x\in\FF_q^s$. The map $y\mapsto x\cdot y$ is onto $\FF_q$, so
  $(q-1)q^{s-1}$ vectors have nonzero image. Dividing by the $q-1$ nonzero representatives of each
  projective point gives the row weight $q^{s-1}$.

  If $\bar{x}\neq\bar{x'}$, then $x$ and $x'$ are linearly independent, so
  $y\mapsto(x\cdot y,x'\cdot y)$ is onto $\FF_q^2$, with fibers of size $q^{s-2}$. There are
  $(q-1)^2$ pairs of nonzero values, giving $(q-1)^2q^{s-2}$ vectors and hence
  $(q-1)q^{s-2}$ projective points.
\end{proof}

For instance, the rows of $N(3,3)$ are the complements of the $13$ lines of the projective plane
$\PG(2,3)$. Each consists of $9$ points and any two share $6$ points, so these $13$ sets form a
$6$-Oddtown on $13$ ground elements.

The preceding lemma controls the Gram entries of $N(s,q)$. We now establish the complementary low-rank
property when $q$ is a power of $p$.

\begin{proposition}\label{pro:projectiveRank}
  Let $p$ be a prime, let $e,\alpha\geq1$, and put $q:=p^e$. There is a constant
  $C_0=C_0(p,e,\alpha)$
  such that
  \[
    \rank_{\ZZ/p^{\alpha}\ZZ}(N(s,q))\leq C_0 s^{e\varphi(p^{\alpha})}
    \qquad\text{for every } s\geq2.
  \]
\end{proposition}

\begin{proof}
  Fix an $\FF_p$-basis of $\FF_q$ and use it to identify $\FF_q^s$ with $\FF_p^{es}$. For every point
  $\bar{x}\in\PG(s-1,q)$ choose a representative $x\in\FF_q^s\setminus\{0\}$. Replace each entry of
  the $\FF_p$-coordinate vector of $x$ by its unique representative in $\{0,\dotsc,p-1\}$, and denote
  the resulting vector in $\ZZ^{es}$ by $\widetilde{x}$.

  Let $X=(X_1,\dotsc,X_{es})$ and $Y=(Y_1,\dotsc,Y_{es})$ be disjoint tuples of indeterminates. The map
  $(x,y)\mapsto x\cdot y$ from $\FF_q^s\times\FF_q^s$ to $\FF_q$ is $\FF_p$-bilinear, so in the
  chosen basis its coordinates are bilinear forms
  \[
    L_1(X,Y),\dotsc,L_e(X,Y)\in\FF_p[X,Y].
  \]
  In the chosen basis, the norm map is represented by a homogeneous polynomial
  $\nu\in\FF_p[Z_1,\dotsc,Z_e]$ of degree $e$. Namely, if $(z_1,\dotsc,z_e)$ is the coordinate vector
  of $z\in\FF_q$, then
  \[
    \nu(z_1,\dotsc,z_e)=\operatorname{N}_{\FF_q/\FF_p}(z).
  \]
  This polynomial is the determinant of the matrix of the $\FF_p$-linear map $w\mapsto zw$, whose
  entries are linear forms in $z_1,\dotsc,z_e$. Multiplication by $z$ is invertible exactly when
  $z\neq0$, so $\nu(z_1,\dotsc,z_e)=0$ if and only if $z=0$. Hence
  \[
    P(X,Y):=\nu(L_1(X,Y),\dotsc,L_e(X,Y))\in\FF_p[X,Y]
  \]
  satisfies $\deg_X P\leq e$ and
  \[
    P(x,y)=0
    \qquad\Longleftrightarrow\qquad
    x\cdot y=0.
  \]
  Here and below, polynomial evaluation at $x,y\in\FF_q^s$ means evaluation at their coordinate
  vectors in $\FF_p^{es}$ under the fixed identification.

  Let $\widetilde{P}\in\ZZ[X,Y]$ be obtained by replacing every coefficient of $P$ by its unique
  representative in $\{0,\dotsc,p-1\}$, and put
  \[
    F(X,Y):=\widetilde{P}(X,Y)^{\varphi(p^{\alpha})}.
  \]
  For the chosen representatives $x,y$, we have
  \[
    \widetilde{P}(\widetilde{x},\widetilde{y})\equiv0\pmod p
    \qquad\Longleftrightarrow\qquad
    x\cdot y=0.
  \]
  If $x\cdot y\neq0$, then $p\nmid\widetilde{P}(\widetilde{x},\widetilde{y})$, and Euler's theorem gives
  $F(\widetilde{x},\widetilde{y})\equiv1\pmod{p^{\alpha}}$. If $x\cdot y=0$, then
  $p\mid\widetilde{P}(\widetilde{x},\widetilde{y})$. Since
  $\varphi(p^{\alpha})=p^{\alpha-1}(p-1)\geq\alpha$, it follows that
  $p^{\alpha}\mid F(\widetilde{x},\widetilde{y})$. Therefore
  \begin{equation}\label{eq:projectiveEvaluation}
    F(\widetilde{x},\widetilde{y})\equiv N(s,q)_{\bar{x},\bar{y}}\pmod{p^{\alpha}}.
  \end{equation}

  Put $D:=e\varphi(p^{\alpha})$. Since $\deg_X F\leq D$, there are polynomials
  $F_{\beta}\in\ZZ[Y]$ such that
  \[
    F(X,Y)=\sum_{\substack{\beta\in\ZZ_{\geq0}^{es}\\\abs{\beta}\leq D}}
      X^{\beta}F_{\beta}(Y).
  \]
  Let $\mathcal{M}:=\{\beta\in\ZZ_{\geq0}^{es}:\abs{\beta}\leq D\}$, and write
  $\widetilde{x}^{\beta}:=\prod_{r=1}^{es}\widetilde{x}_r^{\beta_r}$. Define matrices $A$ and $B$
  over $\ZZ/p^{\alpha}\ZZ$ by reducing the following integer entries modulo $p^{\alpha}$:
  \[
    A_{\bar{x},\beta}:=\widetilde{x}^{\beta},
    \qquad
    B_{\beta,\bar{y}}:=F_{\beta}(\widetilde{y})
    \qquad(\beta\in\mathcal{M}).
  \]
  The polynomial expansion and \eqref{eq:projectiveEvaluation} give $N(s,q)=AB$ over
  $\ZZ/p^{\alpha}\ZZ$. Therefore
  \[
    \rank_{\ZZ/p^{\alpha}\ZZ}(N(s,q))\leq\abs{\mathcal{M}}
      =\binom{es+D}{D}
      \leq(es+D)^D
      \leq(e+D)^D s^D.
  \]
  Taking $C_0:=(e+D)^D$, where $D=e\varphi(p^{\alpha})$, proves the claimed bound.
\end{proof}

\subsection{Blocks of a prescribed order}\label{sec:blocks}
We now combine the projective matrices so as to distinguish one prime-power factor of $\ell$. Fix the
factorization $\ell=p_1^{\alpha_1}\dotsb p_{\omega}^{\alpha_{\omega}}$, where $\omega\geq2$. For
$i\neq j$, let $e_{ij}$ be the multiplicative order of $p_j$ modulo $p_i^{\alpha_i}$, and put
\begin{equation}\label{eq:qij}
  q_{ij}:= p_j^{e_{ij}},
  \qquad\text{so that}\qquad
  q_{ij}\equiv1 \pmod{p_i^{\alpha_i}}.
\end{equation}
For fixed $i$, we take one factor $N(s_j,q_{ij})$ for every $j\neq i$, with $s_j$ chosen large enough
for the required powers of $p_j$ to occur. Since $q_{ij}$ is a power of $p_j$, the $j$-th factor then
supplies divisibility by $p_j^{\alpha_j}$. On the other hand,
$q_{ij}\equiv1\pmod{p_i^{\alpha_i}}$, so the tensor rows have weight $1$ modulo
$p_i^{\alpha_i}$, while a factor $q_{ij}-1$ supplies divisibility by $p_i^{\alpha_i}$ whenever two
tensor rows differ. Moreover, modulo $p_t^{\alpha_t}$ with $t\neq i$, the $t$-th factor has small
inner rank.

\begin{proposition}\label{pro:tensorBlock}
  There are constants $c,C_1>0$ and $R_0\geq2$, depending only on $\ell$, such that for every
  $i\in[\omega]$ and every real $R\geq R_0$ there is a square $0$--$1$ matrix $B$ of order $d$ whose
  rows form an $\ell$-Oddtown, with
  \begin{equation}\label{eq:tensorOrder}
    cR\leq d\leq R,
  \end{equation}
  and such that, for every $t\neq i$,
  \begin{equation}\label{eq:tensorRank}
    \rank_{\ZZ/p_t^{\alpha_t}\ZZ}(B)
    \leq C_1 R^{\frac{\omega-2}{\omega-1}}
      (\log R)^{e_{it}\varphi(p_t^{\alpha_t})}.
  \end{equation}
\end{proposition}

\begin{proof}
  Fix $i\in[\omega]$ and $R\geq R_0$. For each $j\neq i$, choose an integer $s_j\geq2$ as large as
  possible subject to
  \[
    q_{ij}^{s_j}\leq R^{1/(\omega-1)},
  \]
  and put
  \[
    B:=\bigotimes_{j\neq i}N(s_j,q_{ij}).
  \]
  By enlarging $R_0$ if necessary, we may assume that
  \begin{equation}\label{eq:tensorLowerBound}
    s_j\geq 2+\left\lceil\frac{\alpha_j}{e_{ij}}\right\rceil
    \qquad (j\neq i).
  \end{equation}

  The order of $B$ is
  \[
    d=\prod_{j\neq i}\frac{q_{ij}^{s_j}-1}{q_{ij}-1}.
  \]
  Maximality gives $q_{ij}^{s_j+1}>R^{1/(\omega-1)}$, and hence
  \[
    q_{ij}^{-2}R^{1/(\omega-1)}
    <q_{ij}^{s_j-1}
    \leq\frac{q_{ij}^{s_j}-1}{q_{ij}-1}
    \leq q_{ij}^{s_j}
    \leq R^{1/(\omega-1)}.
  \]
  Multiplying over $j\neq i$ proves \eqref{eq:tensorOrder}, after taking
  \[
    c:=\min_{r\in[\omega]}\prod_{j\neq r}q_{rj}^{-2}>0.
  \]

  Every row of $B$ has the form
  \[
    r=\bigotimes_{j\neq i}r_j,
  \]
  where $r_j$ is a row of $N(s_j,q_{ij})$. By \Cref{lem:projectiveCounts}, the weight of $r$ is
  $\prod_{j\neq i}q_{ij}^{s_j-1}$. For $j\neq i$, the factor
  $q_{ij}^{s_j-1}=p_j^{e_{ij}(s_j-1)}$ is divisible by $p_j^{\alpha_j}$ by
  \eqref{eq:tensorLowerBound}, while $q_{ij}\equiv1\pmod{p_i^{\alpha_i}}$. Thus the row weight is
  congruent to $1$ modulo $p_i^{\alpha_i}$, and in particular is not divisible by $\ell$.

  Now let
  \[
    r:=\bigotimes_{j\neq i}r_j,
    \qquad
    r':=\bigotimes_{j\neq i}r'_j
  \]
  be two distinct rows of $B$.
  Their dot product factors as
  \[
    r\cdot r'=\prod_{j\neq i}(r_j\cdot r'_j).
  \]
  By \Cref{lem:projectiveCounts},
  \[
    r_j\cdot r'_j=
    \begin{cases}
      q_{ij}^{s_j-1},&r_j=r'_j,\\
      (q_{ij}-1)q_{ij}^{s_j-2},&r_j\neq r'_j.
    \end{cases}
  \]
  Both values are divisible by $p_j^{\alpha_j}$ by \eqref{eq:tensorLowerBound}. Since $r\neq r'$,
  we have $r_{j_0}\neq r'_{j_0}$ for some $j_0\neq i$. The corresponding factor is divisible
  by $q_{ij_0}-1$, and hence by $p_i^{\alpha_i}$ by \eqref{eq:qij}. Thus $r\cdot r'$ is divisible by
  every prime-power factor of $\ell$, and therefore by $\ell$.

  Finally fix $t\neq i$. Since $q_{it}=p_t^{e_{it}}$, \Cref{pro:projectiveRank} gives
  \[
    \rank_{\ZZ/p_t^{\alpha_t}\ZZ}\bigl(N(s_t,q_{it})\bigr)
      \leq C_0 s_t^{e_{it}\varphi(p_t^{\alpha_t})}.
  \]
  Every other tensor factor has inner rank at most its order, which is at most
  $R^{1/(\omega-1)}$. There are $\omega-2$ such factors. Since $s_t=O_{\ell}(\log R)$,
  \Cref{pro:rankProperties}(iii) gives \eqref{eq:tensorRank}, after choosing $C_1=C_1(\ell)$ uniformly
  over the finitely many pairs $(i,t)$.
\end{proof}

The preceding proposition gives a block whose order is comparable to a prescribed scale $R$, but not
necessarily equal to it. To obtain every sufficiently large order $m$, we place successive blocks on the
diagonal until only a bounded remainder remains. Block-diagonal sums preserve the Oddtown condition,
and their inner ranks are at most the sums of the block ranks. The remainder decreases geometrically,
so only $O_{\ell}(\log m)$ blocks are needed.

\begin{proposition}\label{pro:prescribedOrder}
  There are constants $n_0,C_2,C>0$, depending only on $\ell$, such that, for every $i\in[\omega]$ and
  every integer $m\geq n_0$, there is an $m$-by-$m$ $0$--$1$ matrix $M_i$ whose rows form an
  $\ell$-Oddtown and such that
  \begin{equation}\label{eq:blockRank}
    \rank_{\ZZ/p_t^{\alpha_t}\ZZ}(M_i)
    \leq C_2m^{\frac{\omega-2}{\omega-1}}(\log m)^{C}
    \qquad(t\neq i).
  \end{equation}
\end{proposition}

\begin{proof}
  Set $R_1:=m$. As long as $R_a\geq R_0$, apply \Cref{pro:tensorBlock} to obtain a matrix
  $B_a$ of order $d_a$, and set $R_{a+1}:=R_a-d_a$. Since $R_{a+1}\leq(1-c)R_a$, the process stops
  after $O_{\ell}(\log m)$ steps with $R_{K+1}<R_0$. Let $M_i$ be the block-diagonal matrix with
  blocks $B_1,\dotsc,B_K$, followed by $R_{K+1}$ copies of the $1$-by-$1$ identity matrix. Its rows
  form an $\ell$-Oddtown because the blocks have disjoint supports, the rows of each $B_a$ form an
  $\ell$-Oddtown, and every remaining row has weight $1$.

  Fix $t\neq i$. For the block extracted at remainder $R_a$, \eqref{eq:tensorRank} gives
  \[
    \rank_{\ZZ/p_t^{\alpha_t}\ZZ}(B_a)
      =O_{\ell}\left(
        (\log m)^{e_{it}\varphi(p_t^{\alpha_t})}
        R_a^{\frac{\omega-2}{\omega-1}}
      \right).
  \]
  Summing over the $O_{\ell}(\log m)$ blocks, using $R_a\leq m$, and adding the bounded
  contribution of the identity blocks proves \eqref{eq:blockRank}, with
  \[
    C:=1+\max_{\substack{r,t\in[\omega]\\r\neq t}}e_{rt}\varphi(p_t^{\alpha_t}).
  \]
  Increasing $C_2$ if necessary makes the bound uniform in $i$ and $t$.
\end{proof}

\subsection{Superimposing the blocks}\label{sec:superimpose}
For each $i$, \Cref{pro:prescribedOrder} provides an $\ell$-Oddtown matrix $M_i$ of the same order $m$.
Placing these matrices on the same $m$ columns preserves their diagonal Gram blocks, but leaves the
mixed blocks $M_iM_j^{\top}$ uncontrolled. We will append a small number of columns that cancel all
these mixed blocks modulo $\ell$. For each prime-power factor of $\ell$, at least one of $M_i,M_j$ has
small inner rank, and hence so does $M_iM_j^{\top}$. The Chinese remainder theorem first combines
these separate bounds into a rank bound modulo $\ell$.

\begin{lemma}\label{lem:crtRank}
  For every integer matrix $W$,
  \[
    \rank_{\ZZ/\ell\ZZ}(W)
      =\max_{t\in[\omega]}\rank_{\ZZ/p_t^{\alpha_t}\ZZ}(W).
  \]
\end{lemma}

\begin{proof}
  Reducing a factorization over $\ZZ/\ell\ZZ$ modulo $p_t^{\alpha_t}$ gives
  \[
    \rank_{\ZZ/p_t^{\alpha_t}\ZZ}(W)
      \leq\rank_{\ZZ/\ell\ZZ}(W)
    \qquad(t\in[\omega]),
  \]
  and hence
  \[
    \max_{t\in[\omega]}\rank_{\ZZ/p_t^{\alpha_t}\ZZ}(W)
      \leq\rank_{\ZZ/\ell\ZZ}(W).
  \]
  Conversely, for each $t\in[\omega]$, choose a factorization
  $W=X_tY_t$ over $\ZZ/p_t^{\alpha_t}\ZZ$ of inner dimension
  $\rank_{\ZZ/p_t^{\alpha_t}\ZZ}(W)$. Add zero columns to $X_t$ and zero rows to $Y_t$ so that all
  these factorizations have width
  \[
    \max_{t\in[\omega]}\rank_{\ZZ/p_t^{\alpha_t}\ZZ}(W).
  \]
  The Chinese remainder theorem, applied entrywise, gives matrices $X_0$ and $Y_0$ over
  $\ZZ/\ell\ZZ$ whose reductions modulo $p_t^{\alpha_t}$ are $X_t$ and $Y_t$, respectively, for
  every $t$. Thus
  $X_0Y_0$ and $W$ have the same reduction modulo every prime-power factor of $\ell$, so
  $X_0Y_0=W$ over $\ZZ/\ell\ZZ$. Therefore
  \[
    \rank_{\ZZ/\ell\ZZ}(W)
      \leq\max_{t\in[\omega]}\rank_{\ZZ/p_t^{\alpha_t}\ZZ}(W),
  \]
  which is the reverse inequality.
\end{proof}

The next lemma turns a low-rank integer matrix $W$ into binary correction columns: it produces matrices
$U$ and $V$ whose diagonal Gram matrices vanish modulo $\ell$ and whose mixed Gram matrix is $-W$.

\begin{lemma}\label{lem:gadget}
  Let $W$ be an integer matrix with $a$ rows and $b$ columns, and put
  \[
    r:=\rank_{\ZZ/\ell\ZZ}(W).
  \]
  There are $0$--$1$ matrices $U$ with $a$ rows and $V$ with $b$ rows, having the same number of
  columns, such that
  \begin{equation}\label{eq:gadget}
    UU^{\top}\equiv0\pmod{\ell},
    \qquad
    VV^{\top}\equiv0\pmod{\ell},
    \qquad
    UV^{\top}\equiv-W\pmod{\ell},
  \end{equation}
  and their common number of columns is at most $(2\ell-1)(\ell-1)^{2}r$.
\end{lemma}

\begin{proof}
  If $r=0$, take $U$ and $V$ to have no columns. Hence assume $r\geq1$. Since $-1$ is a unit modulo
  $\ell$, the matrices $W$ and $-W$ have the same inner rank, so choose a factorization
  \[
    -W=XY
    \qquad\text{over }\ZZ/\ell\ZZ,
  \]
  where $X$ has $r$ columns and $Y$ has $r$ rows. For each $h\in[r]$, let $x_h$ be the $h$-th column
  of $X$, and let $y_h$ be the transpose of the $h$-th row of $Y$. Choose representatives in
  $\{0,\dotsc,\ell-1\}$ for their coordinates and regard the resulting vectors as integer vectors.
  Each $x_h$ and $y_h$ is a sum of at most $\ell-1$ $0$--$1$ vectors: for each
  $j=1,\dotsc,\ell-1$, take the indicator of the coordinates whose representatives are at least $j$.
  Reducing these decompositions modulo $\ell$ and expanding
  $-W=\sum_{h=1}^{r}x_hy_h^{\top}$, we see that, over $\ZZ/\ell\ZZ$, the matrix $-W$ is a sum of at
  most $r(\ell-1)^2$ matrices $uv^{\top}$ with
  $u\in\{0,1\}^{a}$ and $v\in\{0,1\}^{b}$.

  Construct $U$ and $V$ by concatenating, for each such summand $uv^{\top}$, one copy of the column
  pair $(u,v)$, then $\ell-1$ copies of $(u,0)$, and finally $\ell-1$ copies of $(0,v)$. These
  $2\ell-1$ columns contribute
  \[
    \ell uu^{\top},
    \qquad
    \ell vv^{\top},
    \qquad
    uv^{\top}
  \]
  to $UU^{\top}$, $VV^{\top}$, and $UV^{\top}$, respectively. Summing over all the outer products
  proves \eqref{eq:gadget}, and uses at most $(2\ell-1)(\ell-1)^2r$ columns.
\end{proof}

To finish, we choose $m$ just below $n$, apply the gadget to every mixed block, and place its two parts
in the corresponding pair of row blocks.

\begin{proof}[Proof of \Cref{thm:construction}]
  Put
  \[
    C_3:=\binom{\omega}{2}(2\ell-1)(\ell-1)^2C_2.
  \]
  Let $n$ be sufficiently large in terms of $\ell$, and set
  \[
    m:=n-\Bigl\lceil C_3n^{\frac{\omega-2}{\omega-1}}(\log n)^{C}\Bigr\rceil.
  \]
  Since $n^{-1/(\omega-1)}(\log n)^{C}=o(1)$, we have $m\geq n/2\geq n_0$. Choose matrices
  $M_1,\dotsc,M_{\omega}$ of order $m$ as in \Cref{pro:prescribedOrder}, and define the large
  $\omega m$-by-$m$ matrix
  \[
    M:=\begin{pmatrix}M_1\\ \vdots\\ M_{\omega}\end{pmatrix}.
  \]
  Define an $\omega m$-by-$\omega m$ block matrix $W=(W_{ij})_{i,j\in[\omega]}$ by
  \[
    W_{ij}:=
    \begin{cases}
      M_iM_j^{\top},&i<j,\\
      0,&i\geq j.
    \end{cases}
  \]
  Thus
  \begin{equation}\label{eq:largeW}
    MM^{\top}
      =\diag\bigl(M_1M_1^{\top},\dotsc,M_{\omega}M_{\omega}^{\top}\bigr)+W+W^{\top}.
  \end{equation}

  Fix $i<j$. For every $t\in[\omega]$, \Cref{pro:rankProperties}(i) and (ii) give
  \[
    \rank_{\ZZ/p_t^{\alpha_t}\ZZ}(W_{ij})
      \leq\min\left\{
        \rank_{\ZZ/p_t^{\alpha_t}\ZZ}(M_i),
        \rank_{\ZZ/p_t^{\alpha_t}\ZZ}(M_j)
      \right\}.
  \]
  At least one of $i,j$ differs from $t$, so \eqref{eq:blockRank} bounds the right-hand side by
  \[
    C_2m^{\frac{\omega-2}{\omega-1}}(\log m)^{C}.
  \]
  By \Cref{lem:crtRank}, the same bound holds for $\rank_{\ZZ/\ell\ZZ}(W_{ij})$. Apply
  \Cref{lem:gadget} to obtain $0$--$1$ matrices $U_{ij}$ and $V_{ij}$, with a common number $s_{ij}$ of
  columns, such that
  \[
    U_{ij}U_{ij}^{\top}\equiv V_{ij}V_{ij}^{\top}\equiv0,
    \qquad
    U_{ij}V_{ij}^{\top}\equiv-W_{ij}
    \pmod{\ell}.
  \]
  Summing the column bounds in \Cref{lem:gadget} over the $\binom{\omega}{2}$ pairs gives
  \begin{equation}\label{eq:totalAppended}
    s:=\sum_{i<j}s_{ij}
      \leq C_3m^{\frac{\omega-2}{\omega-1}}(\log m)^{C}
      \leq C_3n^{\frac{\omega-2}{\omega-1}}(\log n)^{C}
      \leq n-m.
  \end{equation}

  Define the large $\omega m$-by-$s$ matrices $U$ and $V$ with row blocks indexed by $k\in[\omega]$
  and column blocks indexed by the pairs $i<j$ as follows:
  \[
    U_{k,(i,j)}:=
    \begin{cases}
      U_{ij},&k=i,\\
      0,&k\neq i,
    \end{cases}
    \qquad
    V_{k,(i,j)}:=
    \begin{cases}
      V_{ij},&k=j,\\
      0,&k\neq j.
    \end{cases}
  \]
  The nonzero blocks of $U$ and $V$ lie in different row blocks within each column block, so $U+V$ is
  a $0$--$1$ matrix. Since distinct pairs occupy disjoint column blocks, the gadget congruences give
  \begin{equation}\label{eq:largeGadget}
    UU^{\top}\equiv VV^{\top}\equiv0,
    \qquad
    UV^{\top}\equiv-W
    \pmod{\ell}.
  \end{equation}

  By \eqref{eq:totalAppended}, the matrix
  \[
    \widehat M:=\bigl(M\mid U+V\mid 0_{\omega m\times(n-m-s)}\bigr)
  \]
  is an $\omega m$-by-$n$ $0$--$1$ matrix. By \eqref{eq:largeW} and \eqref{eq:largeGadget},
  \begin{align*}
    \widehat M\widehat M^{\top}
      &\equiv MM^{\top}+(U+V)(U+V)^{\top}\\
      &\equiv MM^{\top}-W-W^{\top}\\
      &=\diag\bigl(M_1M_1^{\top},\dotsc,M_{\omega}M_{\omega}^{\top}\bigr)
      \pmod{\ell}.
  \end{align*}
  Since the rows of every $M_i$ form an $\ell$-Oddtown, the displayed congruence shows that
  $\widehat M$ is an $\ell$-Oddtown matrix with
  \[
    \omega m\geq \omega n-O_{\ell}\Bigl(n^{\frac{\omega-2}{\omega-1}}(\log n)^{C}\Bigr)
  \]
  members on $n$ ground elements. This proves the theorem.
\end{proof}

\begin{remark}\label{rem:omegaTwo}
  When $\omega=2$, the exponent $\frac{\omega-2}{\omega-1}$ vanishes, so \Cref{thm:construction} gives
  $f_{\ell}(n)\geq2n-O_{\ell}((\log n)^{C})$. For $\ell=6$, taking $p_1=2$ and $p_2=3$, the exponent
  $C$ in \eqref{eq:blockRank} may be taken to be $3$: the maximum term is $2$, since
  $e_{12}\varphi(3)=e_{21}\varphi(2)=2$, and the additional $1$ comes from the
  $O_{\ell}(\log m)$ block-diagonal summation. Together with \Cref{thm:main}, this gives
  \[
    2n-O\bigl((\log n)^{3}\bigr)\leq f_6(n)\leq 2n-4\log n+11.
  \]
\end{remark}

\section{The upper bounds}\label{sec:upper}

\subsection{Algebraic part: proof of the \texorpdfstring{$2\omega \log n$}{2 omega log n} bound}\label{sec:algebra}
The following lemma is easy for prime moduli using the usual linear algebra in the prime field $\FF_p$.
The fact that it holds for any prime power is the key to proving \Cref{thm:main} in full generality.
\begin{lemma}\label{lem:dimension}
  Suppose that $u_1,\ldots,u_{n-k}\in \ZZ^n$ and $v_1,\dotsc,v_m\in \ZZ^n$ are vectors satisfying
  \begin{equation}\label{eq:orthoConditions}
  \begin{alignedat}{2}
    u_i\cdot u_i&\nequiv 0\pmod{p^{\alpha}}&\qquad&\text{for all }i,\\
    u_i\cdot u_j&\equiv 0\pmod{p^{\alpha}}&&\text{for all }i\neq j,\\
    v_i\cdot v_j&\equiv 0\pmod{p^{\alpha}}&&\text{for all }i,j,\\
    u_i\cdot v_j&\equiv 0\pmod{p^{\alpha}}&&\text{for all }i,j.
  \end{alignedat}
  \end{equation}
  Then the vector space spanned by $v_1\bmod{p},\ldots,v_m\bmod{p}$ over $\FF_p$ has dimension at most $k/2$.
\end{lemma}

We postpone the proof of this lemma until the end of this subsection. For now, we show how to use the lemma to prove
that every $\ell$-Oddtown $\cA\subset 2^{[n]}$ satisfies
$\abs{\cA}\leq \omega n-2\omega \log n+11$.

We need the following well-known fact (appearing, for example, in \cite{odlyzko1981ranks}).
\begin{lemma}\label{lem:schwartzZippel}
  Let $\FF$ be a field, and $V\subset \FF^n$ be a subspace of dimension $d$. Then
  $\abs{V\cap \{0,1\}^n}\leq 2^d$.
\end{lemma}
\begin{proof}
  Let $M$ be a $d$-by-$n$ matrix of rank $d$ whose rows span $V$.
  Since the row rank is equal to the column rank, some set of $d$ columns, say the first $d$ columns,
  span the column space of $M$. This means that every point in $V$ is determined by its first $d$
  coordinates, which implies the lemma.
\end{proof}

Let $\ell=p_1^{\alpha_1}\dotsb p_{\omega}^{\alpha_\omega}$ be the prime factorization of $\ell$.
Recall that $\cA_i:=\{A\in \cA: \abs{A}\nequiv 0\pmod{p_i^{\alpha_i}}\}$, and set $\cA_i':=\cA\setminus \cA_i$.

Fix any $i\in [\omega]$, and apply \Cref{lem:dimension} with the characteristic functions
of sets in $\cA_i$ and in $\cA_i'$ playing the role of vectors $u_1,\dotsc,u_{n-k}$ and $v_1,\dotsc,v_m$ respectively.
We learn that
\begin{equation}\label{eq:rankSizeBound}
  \dim(\spn{p_i} \cA_i')\leq (n-\abs{\cA_i})/2.
\end{equation}
  In view of \Cref{lem:schwartzZippel}, this implies that
\begin{equation}\label{eq:purealg}
  \abs{\cA_i'}\leq 2^{(n-\abs{\cA_i})/2}.
\end{equation}

If $\omega\geq 3$ and $\abs{\cA_i'}\leq n$ for some $i$, then $\abs{\cA}\leq \abs{\cA_i}+\abs{\cA_i'}\leq 2n\leq \omega n - 2\omega\log n+11$.
If $\omega\geq 3$ and $\abs{\cA_i'}\geq n$ for all $i$, then \eqref{eq:purealg}
implies that $\abs{\cA}\leq \sum_i \abs{\cA_i}\leq \sum_i (n- 2\log  n)=\omega n - 2\omega\log n$.

Similarly, if $\omega=2$ and $\abs{\cA_i'}\leq n/2$ for some $i=1,2$, then $\abs{\cA}\leq \tfrac{3}{2}n \leq 2n-4\log n+11$. Otherwise, for both $i=1,2$ we have $\abs{\cA_i}\leq n-2\log \abs{\cA'_i}\leq n-2\log n+2$,
implying $\abs{\cA}\leq 2n-4\log n+4$ in this case. So, in either case $\abs{\cA}\leq 2n-4\log n+11$ holds.

\begin{proof}[Proof of \Cref{lem:dimension}.]
    We first observe that the properties of the $u_i$ and $v_j$ vectors are preserved under the following operations:
    \begin{enumerate}[(i)]
        \item Replacing any vector by its sum or difference with a $v_j$; \label{ope:one}
        \item Multiplying any vector by an integer coprime to $p$; \label{ope:two}
        \item Adding a multiple $p^{\alpha}w$ ($w\in\ZZ^n$) to any vector. \label{ope:three}
    \end{enumerate}

    Let $V:=\spn{p} \{v_1\bmod p,\dotsc,v_m\bmod p\}$.
    By removing some of the vectors, we may assume without loss of generality that $v_1\bmod p,\dotsc,v_m\bmod p$ form a basis for $V$, so that $m=\dim_p V$.
    By Gaussian elimination applied to these $m$ vectors, there are $m$ coordinates onto which $V$ projects bijectively. We may assume that these are the first $m$ coordinates. Using operation (\ref{ope:one}), we can also assume that for all $s,t\in[m]$,\linebreak $v_{s}^{(t)}\not\equiv 0\pmod{p}$ if and only if $s=t$, where $v_s^{(t)}$ denotes the $t$-th coordinate of~$v_s$. Moreover, note that for any integer $z\nequiv 0\pmod{p}$, there exists $z^\dagger$ such that $zz^\dagger\equiv 1\pmod{p^\alpha}$. Thus, we may use operation (\ref{ope:two}) to first make $v_s^{(s)}\equiv 1\pmod{p^\alpha}$, and then use operations (\ref{ope:one}) and (\ref{ope:three}) to ensure that $v_s^{(t)}=\delta_{s,t}$.
    
    We may then apply operation (\ref{ope:one}) to zero out $u_i^{(s)}$, for any $i\in [n-k]$ and any $s\in [m]$, by subtracting a suitable multiple of~$v_s$. Hence, we may assume that $u_i^{(s)}=0$ for all $i\in[n-k]$ and all $s\in[m]$.

    Let $u'_i:=u_i|_{[n]\setminus[m]}$ and $v'_j:=v_j|_{[n]\setminus [m]}$ denote the projections of the vectors $u_i$ and $v_j$ onto the remaining $n-m$ coordinates.
    The assumptions \eqref{eq:orthoConditions}
    imply that
      \begin{alignat*}{2}
         u_i'\cdot u_i'&\nequiv \makebox[\widthof{$-1$}][l]{0} \pmod{p^\alpha} &\quad& \text{for all }i,\\
        u_i'\cdot u_j'&\equiv \makebox[\widthof{$-1$}][l]{0}\pmod{p^\alpha}&\quad& \text{for all }i\neq j,\\
        u_i'\cdot v_j'&\equiv \makebox[\widthof{$-1$}][l]{0}\pmod{p^\alpha}&\quad& \text{for all }i, j,\\
        v_i'\cdot v_i'&\equiv -1\pmod{p^\alpha}&\quad& \text{for all }i,\\
        v_i'\cdot v_j'&\equiv \makebox[\widthof{$-1$}][l]{0} \pmod{p^\alpha}&\quad& \text{for all }i\neq j.
        \end{alignat*}
    
    Supposing on the contrary that $n-m< n+m-k$, we deduce that $\{u_i',v_j'\}_{i\in[n-k],j\in[m]}$ is linearly dependent over $\ZZ$ because they are vectors in $\ZZ^{n-m}$. Let \[
    a_1u_1'+\cdots+a_{n-k}u_{n-k}'=b_1v_1'+\cdots+b_{m}v_m',
  \]
  be a linear dependence, with $a_i,b_j\in \ZZ$. We may assume that $p$ does not divide all $a_i$'s and all $b_j$'s.
  Taking dot product of both sides with $v_j'$ yields $b_j\equiv 0 \pmod{p^\alpha}$, for any $j\in[m]$. Hence,
  \[
    a_1u_1'+\cdots+a_{n-k}u_{n-k}'\equiv 0\pmod{p^\alpha},
  \] is a non-trivial linear combination of the vectors $u_1',\dotsc,u_{n-k}'$ modulo $p^{\alpha}$. Upon taking the dot product with $u_j'$, for any $j$ satisfying $a_j\nequiv 0\pmod p$,
  we reach a contradiction.
\end{proof}

\subsection{Proof of \texorpdfstring{\Cref{thm:main2}}{Theorem 1.3}}\label{sec:combinatorics}
    Let $ \cA $ be an $\ell$-Oddtown on the ground set $[n]$. Recall the definitions, for $i\in[\omega]$,
  \begin{align*}
        \cA_i & :=\{A\in\cA : |A|\nequiv 0\pmod{p_i^{\alpha_i}}\},\\
        \cA_i':=\cA\setminus \cA_i & = \{A\in\cA : |A|\equiv 0\pmod{p_i^{\alpha_i}}\}.
    \end{align*}

Similar to the proof of \Cref{thm:main} in \Cref{sec:algebra}, if $\abs{\cA'_i}\leq n/2$, then 
\[\abs{\cA}=\abs{\cA_i}+\abs{\cA'_i}\leq n+n/2\leq \omega n-(2\omega+\varepsilon)\log n\]
when $n$ is large enough, and we are done in this case. Hence, we may assume that $\abs{\cA'_i}> n/2$ from now on.

Let $M_i$ be the $\abs{\cA_i'}$-by-$n$ matrix whose rows are the characteristic vectors of sets in $\cA_i'$.
The following lemma (whose proof we defer) is, if its conditions are satisfied, a strengthening of \Cref{lem:schwartzZippel}
that we used in the proof of \Cref{thm:main}. We wish to apply this lemma to the matrices $M_i$, when $p_i$ is odd.
\begin{lemma}[Main lemma]\label{lem:01matrix}
    Let $p$ be an odd prime and let $M\in \FF_p^{k\times c}$ be a $0$--$1$ matrix with distinct columns and distinct rows. Set $d=\rank_{p}M$ and assume that $36\log p\leq d\leq 3\log c$. Then we have
        \[d\geq \left(1+\frac{1}{36p^2}\right)(\log k-2).\]
\end{lemma}
However, although the rows of $M_i$ are distinct, the columns might be repeated. Therefore, we shall apply \Cref{lem:01matrix} to the submatrix $M_i'$ of $M_i$ obtained by deleting repeated columns. We say that the odd prime $p_i$ dividing $\ell$ is \emph{good} if $M_i$ has at least $n^{0.4}$ distinct columns, otherwise we say that $p_i$ is \emph{bad}. 

If $p_i$ is good, then $d_i:=\rank_{p_i}(M_i)$ satisfies one of the three possible inequalities: 
\begin{itemize}
\item either $d_i\leq 36\log p_i$,
\item or
  \(
    d_i\geq (1+1/36p_i^2)(\log\abs{\cA'_i}-2),
  \)
\item or $d_i\geq 3\log (n^{0.4})$.
\end{itemize}
In the first case, $\abs{\cA'_i}\leq 2^{d_i}\leq p_i^{36}$ contradicting our assumptions that $\abs{\cA'_i}>n/2$ and that $n$ is large enough.
In the second case,
\begin{equation}\label{eq:diLowerBound}
  d_i\geq \left(1+\frac{1}{36p_i^2}\right)(\log\abs{\cA'_i}-2)\geq \left(1+\frac{1}{36p_i^2}\right)(\log n-3)\geq \left(1+\frac{1}{36p_i^2}\right)\log n-4.
\end{equation}
In the third case, $d_i\geq 1.2\log n\geq (1+1/36p_i^2)\log n$, and so \eqref{eq:diLowerBound} holds in this case, too.

Combining \eqref{eq:diLowerBound} with the inequality $n-\abs{\cA_i}\geq 2d_i$ from \eqref{eq:rankSizeBound} in \Cref{sec:algebra}
leads to
\begin{align*}
  \abs{\cA_i}&\leq n-\left(2+\frac{1}{18p_i^2}\right)\log n+8,&&\text{when }p_i\text{ is good}.\\
\intertext{Even when $p_i$ is not good, we may still use the bound \eqref{eq:purealg} that we proved in \Cref{sec:algebra} and get} 
  \abs{\cA_i}&\leq n-2\log \abs{\cA_i'}\leq n-2\log n+2,&&\text{when }p_i\text{ is bad, or }p_i=2.
\end{align*}
  
Summing the inequalities above over all prime divisors of $\ell$, we get
\[\abs{\cA}\leq \sum_{i=1}^\omega \abs{\cA_i}\leq \omega n-\left(2\omega+\sum_{\text{good }p_i}\frac{1}{18p_i^2}\right)\log n+8\omega.\]

The following lemma shows that most odd prime divisors of $\ell$ are good, otherwise we get a better upper bound.
\begin{lemma}
  Either at most one odd prime divisor of $\ell$ is bad, or
   $\abs{\cA}\leq \omega n-n^{0.2}+\omega\ell$.
\end{lemma}
\begin{proof}
  Suppose that odd primes $p_{j_1},p_{j_2}$ are both bad.
    If $\abs{\cA_{j_1}\cap \cA_{j_2}}\geq n^{0.2}$, then we use the trivial bound $\abs{\cA_i}\leq n$ to deduce that $\abs{\cA}=\abs{\bigcup_{i=1}^{\omega} \cA_i}\leq \sum_{i=1}^{\omega}\abs{\cA_i}-\abs{\cA_{j_1}\cap \cA_{j_2}}\leq \omega n-n^{0.2}$.

    Therefore, we may assume that $\abs{\cA_{j_1}\cap \cA_{j_2}}\leq n^{0.2}$. Since $j_1$ is bad, there are at most $n^{0.4}$ distinct columns in the matrix $M_{j_1}$. Therefore, there is a set of elements $S'\subseteq [n]$ of size at least $n^{0.6}$ such that the columns of $M_{j_1}$ indexed by $S'$ are all the same. Similarly, there is a subset $S\subseteq S'$ of size $n^{0.2}$ such that the columns of $M_{j_2}$ indexed by $S$ are all the same.

    Note that the set family $\cB:=\cA'_{j_1}\cup\cA'_{j_2}$ is also an $\ell$-Oddtown. By construction, all elements of $S$ are contained in exactly the same collection of sets in $\mathcal{B}$. Therefore, we may delete any $\ell$ elements of $S$ from all the sets in $\mathcal{B}$ without changing the $\ell$-Oddtown property. By repeating this until fewer than $\ell$ elements of $S$ remain, we obtain an $\ell$-Oddtown supported on at most $n-\abs{S}+\ell\leq n-n^{0.2}+\ell$ elements. By the trivial bound, we have $\abs{\cA'_{j_1}\cup\cA'_{j_2}}\leq \omega (n-n^{0.2}+\ell)$. Thus,
    \[\abs{\cA}\leq\abs{\cA_{j_1}\cap \cA_{j_2}}+\abs{\cA'_{j_1}\cup\cA'_{j_2}}\leq n^{0.2}+\omega (n-n^{0.2}+\ell)\leq \omega n-n^{0.2}+\omega\ell.\qedhere\]
\end{proof}

By this lemma, we may assume without loss of generality that at most one odd prime divisor of $\ell$ is bad. Therefore, 
\[\sum_{\text{good }p_i}\frac{1}{p_i^2}\geq \sum_{p\in P_\text{odd}^*(\ell)}\frac{1}{p^2},\]
and hence 
\[\abs{\cA}\leq \omega n-(2\omega+\varepsilon)\log n,\]
when $n$ is large enough.

\subsection{Analytic part: Fourier-analytic argument}\label{sec:analysis}
\paragraph{Reduction to a strengthening of \Cref{lem:schwartzZippel}.}
In this part, we prove \Cref{lem:01matrix}. We do this by reducing it to the following result.

    \begin{lemma}[Main lemma, probabilistic form]\label{lem:main}
        Let $p$ be an odd prime and let $c,d$ be positive integers such that $36\log p\leq d\leq 3\log c$. Let $L\in \FF_p^{d\times c}$ be a matrix whose columns are distinct vectors in $\FF_p^d$, and let $v$ be uniformly distributed over $\{0,1\}^d$. Then 
        \[
        \mathbb{P}[L^\top v\in\{0,1\}^c] \leq 3\cdot 2^{-d/36p^2}.
        \] 
    \end{lemma}
    
    \begin{proof}[Proof of \Cref{lem:01matrix} assuming \Cref{lem:main}]
    Let $r_1,\dots, r_k$ be the row vectors of $M$. Since the rank of $M$ is $d$, we may assume without loss of generality that $r_1,\dots,r_d$ span the row space. Furthermore, we may pick $d$ columns such that the vectors $r_1,\dots,r_d$ restricted to these columns are linearly independent. Assume without loss of generality that these columns are the first $d$ columns. Let $r'_1,\dots,r'_d$ be the linear combinations of $r_1,\dots,r_d$ such that the projections of $r'_1,\dots,r'_d$ to the first $d$ columns form the identity matrix. 
    
    Now, let $L\in\FF_p^{d\times c}$ denote the matrix with rows $r'_1,\dots,r'_d$. Consider any row $r_x$ of $M$; since $r_x$ is a linear combination of $r_1,\dots,r_d$, it is also a linear combination of $r'_1,\dots,r'_d$ as well, say $r_x=\sum_{i=1}^d a_ir'_i$. Note that the coefficients $a_1,\dotsc,a_d$ are in $\{0,1\}$ because
    the $i$-th coordinate of $r_x$ is $a_i$, for all $i\in [d]$.
    Therefore, the vector $r_x$ is of the form $v^\top L$ with $v\in\{0,1\}^d$. 

   The first \(d\) selected rows already have distinct columns. Indeed, every other row of \(M\) is a linear combination of these rows, so equality of two columns on the selected rows would imply equality on every row of \(M\). Applying an invertible \(d\times d\) row transformation to the selected rows preserves column distinctness, and therefore \(L\) has distinct columns.

    \Cref{lem:main} says that the number of rows of $M$ satisfies
    \[k\leq 2^d\cdot 3\cdot 2^{-d/36p^2},\]
    and hence 
    \[d\geq \left(1+\frac{1}{36p^2}\right)(\log k-2).\qedhere\]
    \end{proof}

\paragraph{Proof of the main lemma.}
We begin by obtaining bounds on the Fourier transform of the probability distribution induced by the action of distinct linear forms on random $0$--$1$ vectors. These Fourier-analytic estimates will enable us to control the number of $0$--$1$ vectors that lie in subspaces of $\FF_p^c$.

Let $p$ be a prime number. For any function $f\colon \FF_p^c\to\CC$, its \emph{Fourier transform} is defined by \[
\hat{f}\colon\FF_p^c\to \CC, \quad \xi\mapsto\sum_{x\in \FF_p^c} f(x)e^{-2\pi i\langle x,\xi \rangle/p},
\] where the inner product $\langle x,\xi\rangle$ takes values in $(-p/2,p/2]$.

\begin{proposition}\label{pro:fourier}
Let $p$ be a prime number and let $L\in \FF_p^{d\times c}$ be a matrix with distinct columns. Suppose that $v\in\{0,1\}^d$ is chosen uniformly at random and let $f(x)$ denote the probability that $L^\top v=x$. Then, for every $\xi\in \FF_p^c$, we have $\lvert \hat{f}(\xi)\rvert\leq e^{-\frac{\pi^2\sigma(\xi)}{2p^2}}$, where $\sigma(\xi):=\sum_{j=1}^d \langle r_j,\xi \rangle^2$ and $r_j$ denotes the $j$-th row of $L$.
\end{proposition}

\begin{proof}
       Let $f_j(x)$ denote the probability that $v^{(j)}r_j^\top = x$. Since $f_j(0)=f_j(r_j^\top)=1/2$ (or $f_j(r_j^\top)=1$ if $r_j=0$) and $f_j(x)=0$ for $x\notin \{0,r_j^\top\}$, the Fourier transform of $f_j$ is  \begin{align*}
\hat{f_j}(\xi)=\sum_{x\in \FF_p^c} f_j(x)e^{-2\pi i\langle x,\xi \rangle/p}=\frac{1}{2}\left(1+e^{-2\pi i\langle r_j, \xi\rangle/p}\right).
\end{align*} Note that \begin{align*}
f(x)=\sum_{x_1+x_2+\cdots+x_d=x}f_1(x_1)f_2(x_2)\cdots f_d(x_d)=f_1*f_2\cdots * f_d (x),
\end{align*} where $*$ denotes convolution. Consequently, the Fourier transform of $f$ satisfies \begin{align*}
    \hat{f}(\xi)=\Wh{f_1*f_2*\cdots*f_d}(\xi)=\prod_{j=1}^d \hat{f_j}(\xi)=\prod_{j=1}^d\frac{1+e^{-2\pi i\langle r_j, \xi\rangle/p}}{2}.
    \end{align*} Taking absolute values and using the inequality $|(1 + e^{-2\theta i})/{2}| = |\cos(\theta)|\leq e^{-\theta^2/2}$ that is valid for \linebreak$\theta\in (-\pi/2,\pi/2]$, we obtain
    \[
    |\hat{f}(\xi)|
    =\prod_{j=1}^d \left|\cos \left(\frac{\langle r_j,\xi \rangle \pi}{p}\right)\right|
    \leq \prod_{j=1}^d e^{-\langle r_j,\xi \rangle^2\pi^2/2p^2}
=e^{-\pi^2\sigma(\xi)/2p^2}.\qedhere
\]
\end{proof}

We say that a matrix $S \in \FF_p^{d\times c}$ with rows $r_1,r_2,\dotsc,r_d$ is \emph{$\sigma$-admissible} if for any nonzero $\xi\in \FF_p^c$, \[
\sum_{i=1}^d \langle r_i, \xi\rangle^2 > \sigma.
\]
We also write $H(q):=-q\log q-(1-q)\log(1-q)$ for the binary entropy function.

Note that if $S$ is $\sigma$-admissible, then any submatrix $ S'\in\FF_p^{d\times c'}$ is also $\sigma$-admissible since we can set the entries of $\xi$ corresponding to the columns that are not in $S'$ to zero.

\begin{lemma}\label{lem:subMatrix}
Let $p$ be an odd prime number, $\sigma>0$ be a real number, and $L \in \FF_p^{d\times c}$ be a matrix with distinct column vectors. There exists a $\sigma$-admissible submatrix $S\in \FF_p^{d\times c'}$ of $L$ with
\[c'\geq \frac{1}{\log p}\left(\log c-\sigma-(\sigma+d)H\left(\frac{\sigma}{\sigma+d}\right)\right)-1.\]
\end{lemma}
\begin{proof}
Consider all submatrices of $L$ that are obtained from $L$ by deleting some columns.
Let $S\in \FF_p^{d\times c'}$ be a maximal $\sigma$-admissible such submatrix. If there is no such submatrix, we set $c'=0$ and $S$ to be the empty matrix.

Suppose $w$ is a column of $L$ that is not in $S$, and let $S'$ be the matrix obtained by appending $w$ to $S$.
It is not $\sigma$-admissible by the maximality of $S$. This implies the existence of a nonzero vector $\eta \in \FF_p^{c' + 1}$ such that
\begin{equation}\label{eq:inadmissibility}
\sum_{i=1}^d\langle r_i', \eta\rangle ^2 \leq \sigma,
\end{equation}
where $r'_1,\dots,r'_d$ are the row vectors of $S'$.

We claim that $\eta_{c'+1}\neq 0$. Indeed, if $\eta_{c'+1}=0$, then the restriction of $\eta$ to the first $c'$ coordinates
would satisfy $\sigma\geq \sum_{i=1}^d\langle r_i', \eta\rangle ^2=\sum_{i=1}^d\langle r_i, \operatorname{proj}(\eta)\rangle^2$, where $\operatorname{proj}(\eta)$ denotes
the projection onto the first $c'$ coordinates; this would contradict the $\sigma$-admissibility of~$S$.
This observation bounds the number of possible choices for $\eta$ by $(p-1)p^{c'} $. Note that this also holds in the case where $c'=0$ since $\eta$ must be a non-zero vector.

Next, we upper bound the number of all possible $w$. For fixed $\eta$ satisfying $\eta_{c'+1}\neq 0$,
the function $w\mapsto (\langle r_1',\eta\rangle,\dotsc,\langle r_d',\eta\rangle)$ is injective. That is because
the $i$-th coordinate of $w$ is determined by $\langle r_i',\eta\rangle$ and~$r_i$, the latter of which is fixed
because $S$ is fixed. By \eqref{eq:inadmissibility}, the vector $(\langle r_1',\eta\rangle,\dotsc,\langle r_d',\eta\rangle)$ satisfies the
inequality
\begin{equation}\label{eq:squareSum}
  \sum_{i=1}^d x_i^2\leq \sigma.
\end{equation}
  Since the number of solutions to the inequality $\sum_{i=1}^d y_i\leq \sigma$ in nonnegative integers $y_1,\dotsc,y_d$
is $\binom{\lfloor \sigma\rfloor+d}{d}\leq\binom{\sigma+d}{\sigma}$, it follows that the number of vectors $(x_1,x_2,\dotsc,x_d)$ satisfying \eqref{eq:squareSum} is at most $2^{\sigma} \binom{\sigma+d}{\sigma}$. Therefore, the number of column vectors $w$ in $L$ that are not in $S$, which equals $c-c'$, is at most $(p-1)p^{c'} 2^{\sigma} \binom{\sigma+d}{\sigma}$.

Therefore, 
\begin{align*}
    c\leq c'+(p-1)p^{c'}2^{\sigma} \binom{\sigma+d}{\sigma}    \leq p^{c'+1}2^{\sigma+(\sigma+d)H\left(\sigma/(\sigma+d)\right)},
\end{align*}
and hence
\[c'\geq \frac{1}{\log p}\left(\log c-\sigma-(\sigma+d)H\left(\frac{\sigma}{\sigma+d}\right)\right)-1.\qedhere\]
\end{proof}

With these results in hand, we are now ready to complete the proof of our main lemma.
\begin{proof}[Proof of \Cref{lem:main}]
    Upon applying \Cref{lem:subMatrix} with $\sigma=d/36$, we learn that there exists a $\sigma$-admissible submatrix $S\in\FF_p^{d\times c'}$ of $L$, where
    \[c'\geq \frac{1}{\log p}\left(\log c-\sigma-(\sigma+d)H\left(\frac{\sigma}{\sigma+d}\right)\right)-1\geq \frac{1}{\log p}\left(\log c-\frac{d}{4}\right)-1\geq \frac{d}{12\log p}-1.\]
    Here, we used the fact that 
    \(\frac{1}{36}+\frac{37}{36}H\left(\frac{1}{37}\right)\leq \frac{1}{4}\).
    
    Let $T\in \FF_p^{d\times c''}$ be any submatrix of $S$ with $c''$ columns, where 
    \[c'':=\Bigl\lceil\frac{d}{12p^2\log p}\Bigr\rceil\leq \frac{d}{12\log p}-1\leq c'.\]
    Note that $T$ is also $\sigma$-admissible.    
    
    Let $f(x)$ denote the probability mass function of $T^\top v$. By the Fourier inversion formula and the bound on the Fourier transforms from \Cref{pro:fourier}, \begin{align*}
        f(x) & =\mathbb{E}_{\xi\in \FF_p^{c''}} \hat{f}(\xi)e^{2\pi i\langle x, \xi \rangle / p}\\
        & \leq \frac{\hat{f}(0)}{p^{c''}}+\frac{p^{c''}-1}{p^{c''}}e^{-\pi^2\sigma/2p^2}\\
        & \leq \frac{1}{p^{c''}}+e^{-\pi^2\sigma/2p^2}.
    \end{align*}
    
    Therefore, \begin{align*}
        \mathbb{P}[L^\top v \in\{0,1\}^c]  \leq \mathbb{P}[T^\top v\in\{0,1\}^{c''}] = \sum_{x\in\{0,1\}^{c''}}f(x)         \leq 2^{c''}\left(\frac{1}{p^{c''}}+e^{-\pi^2\sigma/2p^2}\right).
    \end{align*}
    We may bound these two terms by 
    \[\left(\frac{2}{p}\right)^{c''}=2^{-c''(\log p-1)}\leq 2^{-d(\log p-1)/12p^2\log p}\leq 2^{-d/36p^2}\]
    and
    \begin{align*}
        2^{c''}e^{-\pi^2\sigma/2p^2}\leq 2^{d/12p^2\log p+1}e^{-\pi^2d/72p^2}\leq 2\cdot2^{d/12p^2\log p}2^{-d/12p^2}\leq 2\cdot2^{-d/36p^2}.
    \end{align*}
    Here, we used the estimates  $(\log p-1)/\log p>1/3$ and $e^{\pi^2/36}\geq 2^{1/6}$. Thus, we may conclude that
    \[\mathbb{P}[L^\top v \in\{0,1\}^c]  \leq 3\cdot2^{-d/36p^2}.\qedhere\]
\end{proof}

\printbibliography
\end{document}